\documentclass[11pt]{amsart}
\usepackage{calc,amssymb,amsthm,amsmath,fullpage}
\RequirePackage[dvipsnames,usenames, table]{xcolor}
\usepackage{hyperref}
\hypersetup{
bookmarks,
bookmarksdepth=3,
bookmarksopen,
bookmarksnumbered,
pdfstartview=FitH,
colorlinks,backref,hyperindex,
linkcolor=Sepia,
anchorcolor=BurntOrange,
citecolor=MidnightBlue,
citecolor=OliveGreen,
filecolor=BlueViolet,
menucolor=Yellow,
urlcolor=OliveGreen
}
\usepackage[all]{xy}
\usepackage{alltt}
\usepackage{multicol}  
\usepackage{xspace}
\usepackage{rotating}
\usepackage{mabliautoref}
\usepackage{colonequals}
\input{kmacros3.sty}
\usepackage{stmaryrd}

\usepackage{verbatim}
\usepackage{enumerate}
\begin{document}
\title{{FastMinors} package for {\it Macaulay2}}
\author{Boyana Martinova}
\email{u1056124@utah.edu}
\address{Department of Mathematics, University of Utah, 155 S 1400 E Room 233, Salt Lake City, UT, 84112}
\author{Marcus Robinson}
\email{mrobinso@reed.edu}
\address{Department of Mathematics, Reed College, 3203 Southeast Woodstock Boulevard Portland, Oregon 97202-8199}
\author{Karl Schwede}
\email{schwede@math.utah.edu}
\address{Department of Mathematics, University of Utah, 155 S 1400 E Room 233, Salt Lake City, UT, 84112}
\author{Yuhui Yao}
\email{weiy@math.uchicago.edu}
\address{Department of Mathematics, University of Chicago, Eckhart Hall, 5734 S University Ave, Chicago IL, 60637, 773 702 7100}

\date{\today}

\begin{abstract}
In this article, we present {\tt FastMinors.m2}, a package in {\it Macaulay2} designed to
 introduce new methods focused on computations in function field linear algebra. 
Some key functionality that our 
package offers includes: finding a submatrix of a given rank in a provided matrix (when present), 
verifying that a ring is regular in codimension n, recursively computing the ideals of minors in a 
matrix, and finding an upper bound of the projective dimension of a module.
\end{abstract}

\thanks{Martinova was supported by a University of Utah Mathematics REU fellowship and by the University of Utah ACCESS program.  
Robinson was supported by NSF RTG grant \#1840190.
Schwede was supported by NSF CAREER grant \#1501102, NSF grant \#1801849 and NSF FRG Grant \#1952522 and a fellowship from the Simons Foundation.
Yao was supported by a University of Utah Mathematics REU fellowship.}

\maketitle

\section{Introduction}

We start with some motivation.  
Suppose that $I = (f_1, \dots, f_m) \subseteq k[x_1, \dots, x_n]$ is a prime ideal.  The corresponding variety $X := V(I)$ is nonsingular if and only if $I$ plus the ideal generated by the minors of size $n-\dim X$ of the Jacobian matrix 
\[
\mathrm{Jac}(X) = \left( \begin{array}{c} {\partial f_i \over \partial x_j}\end{array}\right),
\]
generates the unit ideal.  
Unfortunately, even for relatively small values of $m$ and $n$, the number of such submatrices is prohibitive.  Suppose for instance that $n = 10, m = 15$ and $\dim X = 5$.  Then there are  
\[
  {10 \choose 5} \cdot {15 \choose 5} =  756756
\]
such submatrices.
We cannot reasonably compute all of their determinants.
 This package attempts to fix this in several ways.  
\begin{enumerate}
  \item We try to compute just a portion of the determinants, in a relatively smart way.
  \item We offer some tools for computing determinants that are sometimes faster.
\end{enumerate}
Our techniques have also been applied to the related problem of showing that the singular locus has a certain codimension (for example, checking that a variety is R1 in order to prove normality). 
Of course, computing the singular locus is not the only potential application.  We provide a function for giving a better upper bound on the projective dimension of a non-homogeneous module.  Finally, this package has also been applied in the {\tt RationalMaps} {\it Macaulay2} package.  

We provide the following functions:

\begin{itemize}
  \item{} {\tt getSubmatrixOfRank}, which tries to find a submatrix of a given rank, see \autoref{sec.GetSubmatrixOfRank}.
  \item{} {{\tt isRankAtLeast}}, which uses {\tt getSubmatrixOfRank} to try to find lower bounds for the rank of a matrix, see \autoref{sec.IsRankAtLeast}.
  \item{}  {\tt regularInCodimension}, which tries to verify if an integral domain is regular in codimension n, see \autoref{sec.Rn}.
  \item{}  {\tt projDim}, which tries to find upper bounds for the projective dimension of a non-homogeneous module, see
  \autoref{sec.ProjDim}.
  \item{}  {\tt recursiveMinors}, which computes the ideal of minors of a matrix via a recursive cofactor algorithm, as opposed to the included non-recursive cofactor algorithm, see \autoref{sec.RecursiveMinors}.
\end{itemize}

The latest version of this package is available at:
\begin{center}
  {\tt https://github.com/kschwede/M2/blob/master/M2/Macaulay2/packages/FastMinors.m2}
\end{center}


This paper refers to {\tt FastMinors} version 1.2.2, which we hope will appear in version 1.19 of {\it Macaulay2}.  Earlier versions are also available in the {\it Macaulay2} build tree.

\subsection*{Acknowledgements:}  The authors thank David Eisenbud, Dan Grayson, Elo\'isa Grifo and Zhuang He for valuable conversations and feedback.

\section{Finding interesting submatrices}
\label{sec.FindingInterestingSubmatrices}

A lot of the speedups available in the package come down to finding interesting square submatrices of a given matrix. For example, it is often useful to compute a square submatrix whose determinant has small degree.  The idea is that the determinant of this submatrix will be less likely to vanish.

\subsection{How are the submatrices chosen?}

Consider the following matrix defined over $\mathbb{Q}[x,y]$.  
\[
  \left[\begin{array}{ccc}
    x & xy & 0 \\
    xy^2 & x^6 & 0 \\
    0 & x^2 y^3 & xy^4
  \end{array} \right]
\]
Suppose we want to choose a submatrix of size $2 \times 2$.  
Consider the monomial order {\tt Lex} where $x < y$.  We find, in the matrix, the nonzero element of smallest order.  
In this case, that is $x$.  We choose this element to be a part of our submatrix.  Hence our submatrix will include the first row and column as well.
\[
  \left[\begin{array}{>{\columncolor{red!20}}ccc}
    \rowcolor{red!20}
    \cellcolor{red!40}{\bf x} & xy & 0 \\
    xy^2 & x^6 & 0 \\
    0 & x^2 y^3 & xy^4
  \end{array}\right]
\]
To find the next element, we discard that row and column containing this term.  Now, the next smallest 
element with respect to our monomial order is $xy^4$.
\[
  \left[\begin{array}{cc}
    x^6 & 0 \\
    x^2 y^3 & x y^4
  \end{array}\right]
  \;\;\;\;
  \left[\begin{array}{c>{\columncolor{red!20}}c}    
    x^6 & 0 \\    
    \rowcolor{red!20}
    x^2 y^3 & \cellcolor{red!40} {\bf x y^4}
  \end{array}\right]
\]
Since we are only looking for a $2 \times 2$ submatrix, we stop here. 
We have selected the submatrix with rows $0$ and $2$ and columns $0$ and $2$. 
\[
  \left[\begin{array}{ccc}
    x &  0 \\
    0 & x y^4 \\
  \end{array} \right]
\]
The determinant of that submatrix is $x^2 y ^4$.  This happens to be the smallest $2 \times 2$ minor with respect to the given monomial order (which frequently happens, although it is certainly not always the case).

If we chose a different monomial order, we get a different submatrix, with a different determinant.  For example, 
\begin{itemize}  
  \item{} {\tt Lex}, $x > y$.  We obtain the submatrix with rows $0$ and $1$ and columns $0$ and $1$, whose determinant is $x^7 - x y^3$
  \item{} {\tt GRevLex}, $x < y$.  We obtain the submatrix with rows $0$ and $2$ and columns $0$ and $1$, whose determinant is $x^3 y^3$
\end{itemize}
For any of these strategies, in this package, we randomize the order of the variables before choosing a submatrix.

\subsection{Ways of choosing submatrices}
    \label{subsec.WaysOfChoosingSubmatrices}
    In the end we have the following methods to select submatrices.
    \begin{description}        
        \item[{\tt GRevLexLargest}]  Chose entries which are largest with respect to a random {\tt GRevLex} order.
        \item[{\tt GRevLexSmallest}]  Chose nonzero entries which are smallest with respect to a random {\tt GRevLex} order.
        \item[{\tt GRevLexSmallestTerm}]  Chose  nonzero entries which have the smallest terms with respect to a random {\tt GRevLex} order.
        \item[{\tt LexLargest}]  Chose entries which are largest with respect to a random {\tt Lex} order.
        \item[{\tt LexSmallest}]  Chose nonzero entries which is smallest with respect to a random {\tt Lex} order.
        \item[{\tt LexSmallestTerm}]  Chose nonzero entries which have the smallest terms with respect to a random {\tt GRevLex} order.
        \item[{\tt Points}]  Chose a submatrix whose determinant does not vanish at a random point found on a given ideal.
        \item[{\tt Random}]  Chose a random entries.
        \item[{\tt RandomNonzero}]  Chose random nonzero entries.  
    \end{description}    

    However, from the end user's perspective, normally we are finding several minors, and the strategy will combine several of these methods (one typically does not know which method will work best in a given situation).  For instance, the first minor might be selected by {\tt GRevLexSmallest} and the second minor by {\tt Random}.  How to arrange what method is used (with what probability) is described in \autoref{subsec.TheStrategyOption} below.

    We now describe each of these methods for selecting a submatrix in more detail.  Note we have already described {\tt Lex} and given an initial description of {\tt GRevLex}.  

\subsection{{\tt LexSmallestTerm} and {\tt GRevLexSmallestTerm}}
If we have a matrix whose entries are not monomial, then we could reasonably either pick the 
submatrix of smallest entries with respect to our monomial order: {\tt LexSmallest} or {\tt GRevLexSmallest}.
    
    Alternatively, we can pick the submatrix whose entries have the smallest terms via  {\tt LexSmallestTerm} or {\tt GRevLexSmallestTerm}.  

For example, consider the matrix
\[
  \left[\begin{array}{ccc}
    x^2 + y^2 & 0 & xy + 2x \\
    y^4 - x & 0 & 3 x^5 \\
    x^3 & x^4 y^5  - y^8 & 0
  \end{array} \right]
\]
In this case, if we are choosing the entries with smallest terms, we first replace each entry in the matrix 
with the smallest term.  For example, if we are using {\tt LexSmallestTerm} with $x < y$ we would obtain:
\[
  \left[\begin{array}{ccc}
    x^2 & 0 & 2x \\
    -x & 0 & 3 x^5 \\
    x^3 & x^4 y^5 & 0
  \end{array} \right]
\]
Then we proceed as before.  Notice that if there is a tie, it is broken randomly.  

\begin{remark}
    Different strategies work differently on different examples.  When working with a non-homogeneous matrix, with some entries that have constant terms, those entries will always be chosen first in {\tt LexSmallestTerm} or {\tt GRevLexSmallestTerm}, regardless of the monomial order.  On the other hand, for homogeneous matrices, choosing the smallest term is frequently very effective.
  \end{remark}

\subsection{{\tt GRevLexLargest} and {\tt LexLargest}}

While we can imagine uses for these, in most cases these strategies appear to be worse than random.  Indeed, submatrices picked this way seem likely to already vanish everywhere of interest.  

\subsection{{\tt Points}}

Instead of finding interesting submatrices by inspection, we can alternately find submatrices by trying to find rational points.  In that case, typically we are trying to find a submatrix with full rank on a certain subvariety, defined by an ideal $J$.  We use the package {\tt RandomRationalPoints} \cite{RandomRationalPointsSource} to find a (rational) point $Q$ on $V(J)$ (or a point over some finite extension of our base field).  We then evaluate our entire matrix at that point.  Because we now have a matrix over a field, we use the very fast built-in commands to find pivot rows and columns, and thus find a submatrix of the desired rank.  To use this functionality, use the strategy {\tt Points}

Currently, this functionality only works over a finite field.  In characteristic zero, the {\tt Points} strategy returns random submatrices.

For example, if we are working over $\bF_5[x,y,z]$ with an ideal $I = (z^2 y - x(x-z)(x+z))$, with the matrix 
\[
    M = \left[
        \begin{array}{ccc}
            x^2 & xy & 3y^2 \\
            x^3+y^3 & x^2 + z^2 & y^2 + z^2 \\
            x^2*z & z^2*y & y^2*x \\
        \end{array}
    \right]
\]
Suppose we found the point $(2,0,2)$ on this elliptic curve.  We then evaluate our matrix at that point to obtain:
\[
    M = \left[
        \begin{array}{ccc}
            4 & 0 & 0 \\
            3 & 3 & 4 \\
            3 & 0 & 0 \\
        \end{array}
    \right]
\]
We would then identify a submatrix with nonzero determinant, for instance the top left $2\times 2$ submatrix $\left[
    \begin{array}{cc}
        4 & 0 \\
        3 & 3  \\        
    \end{array}
\right]$
and then return the top determinant of the top $2 \times 2$ submatrix of the original matrix:
\[
    \det \left[
        \begin{array}{ccc}
            x^2 & xy \\
            x^3+y^3 & x^2 + z^2 \\
        \end{array}
    \right] = x^4 + x^2 z^2 - x^4y - xy^4.
\]

Typically, the {\tt Points} strategy will find much better submatrices than the heuristic methods {\tt LexSmallest, GRevLexSmallestTerm} etc.   Thus it will need to consider fewer submatrices and compute fewer determinants.  But to find each submatrix, the algorithm will do a substantial amount of work.  In our experience, using the heuristic methods above provides better performance than {\tt Points} in at least half of example applications.  If the user is using a strategy based on points, we make two recommendations.  
    \begin{itemize}
        \item[(i)]  We recommend the user set the option  {\tt MaxMinors} (the maximum number of minors to be computed) to a relatively low number.
        \item[(ii)]  We recommend setting the option {\tt CodimCheckFunction} to a linear function (such as {\tt t -> t}).  This will force the dimension of the ideal of minors computed so far to be checked more frequently (in our example, after adding every new minor to the ideal of minors)
    \end{itemize}
    For more discussion of {\tt MaxMinors}, {\tt CodimCheckFunction} and other options, see the tutorial {\tt RegularInCodimensionTutorial} in the package documentation.

\subsection{{\tt Random} and {\tt RandomNonzero}}
\begin{description}
  \item[{\tt Random}]  With this strategy, a random submatrix is chosen.
  \item[{\tt RandomNonzero}]  With this strategy, a random nonzero element is chosen in each step following the method used by the other strategies. This guarantees a submatrix where no row or column is zero which can be very useful when dealing with relatively sparse matrices.
\end{description}

\subsubsection{More on {\tt GRevLex}: modifying the underlying matrix}
Finally, when using {\tt GRevLexSmallest} and {\tt GRevLexSmallestTerm} methods, we periodically change the underlying matrix by replacing terms of small order with terms of larger order in order to avoid re-computing the same submatrix.  For example, in the following matrix, after 
    several iterations, we might replace the $x^2$ term with 
\[x^2 \cdot (\text{a random degree 1 polynomial}). \]
It might look something like the following.
\[
  \left[\begin{array}{ccc}
    x^2 & 0 & xy \\
    y^4 & 0 & x^5 \\
    x^3 & x^4 y^5 & 0
  \end{array} \right]
  \rightarrow
  \left[\begin{array}{ccc}
    x^2 (2x-7y) & 0 & xy \\
    y^4 & 0 & x^5 \\
    x^3 & x^4 y^5 & 0
  \end{array} \right].
\]
This forces the algorithm to make different choices.  After several minors are selected, the matrix is reset again to its original form.

\section{{\tt chooseGoodMinors} and controlling how submatrices are selected}  
\label{sec.OptionsForChoosingSubmatrices}

  The function {\tt chooseGoodMinors} tries to choose interesting submatrices of a given matrix.  This is done by running the command:

{{\small\color{blue}
\begin{verbatim}
i1 : R = QQ[x, y, z];
    chooseGoodMinors()
\end{verbatim}}
}
\subsection{The {\tt Strategy} option}
\label{subsec.TheStrategyOption}
The core features included in the package allow the user to choose which 
methods from \autoref{subsec.WaysOfChoosingSubmatrices} should be used when selecting submatrices.  This is done most easily by setting a {\tt Strategy} option to one of the ways of choosing submatrices as above:  {\tt GRevLexSmallest, GRevLexSmallestTerm, GRevLexLargest, LexSmallest, LexSmallestTerm, LexLargest, Points, Random, RandomNonzero}.  However, most of the time it is best to choose several strategies simultaneously, as one doesn't know which strategy will perform the best (in some cases, a combination works best).  Hence instead of choosing a strategy which uses only one method, by default we use several.  Thus you can set the {\tt Strategy} option to one of the following.

\begin{itemize}
  \item {\tt StrategyDefault}: This strategy uses  
  {\tt LexSmallest, LexSmallestTerm, GRevLexSmallest, GRevLexSmallestTerm, Random,} and {\tt RandomNonzero} with equal probability.
  \item {\tt StrategyDefaultNonRandom}: This uses  
  {\tt LexSmallest, LexSmallestTerm, GRevLexSmallest,} and {\tt GRevLexSmallestTerm}  with equal probability.
  \item {\tt StrategyDefaultWithPoints}: This strategy uses {\tt Points} one third of the time and then 
  {\tt LexSmallest, LexSmallestTerm, GRevLexSmallest, GRevLexSmallestTerm} with equal probability the rest of the time.
  \item {\tt StrategyLexSmallest}: chooses 50\% of the submatrices using {\tt LexSmallest} and 50\% 
  using {\tt LexSmallestTerm}.
  \item {\tt StrategyGRevLexSmallest}: chooses 50\% of the submatrices using {\tt GRevLexSmallest} 
  and 50\% using {\tt GRevLexLargest}.
  \item {\tt StrategyPoints}: choose submatrices by finding rational points, evaluating the submatrix at that point, and then doing a computation.
  \item {\tt StrategyRandom}: chooses submatrices by using 50\% {\tt Random} and 50\% {\tt RandonNonzero}.
  \end{itemize}  

The user can also create their own custom strategy.  One creates a {\tt HashTable} which has the following keys 
{\tt LexLargest, LexSmallestTerm, LexSmallest, GRevLexSmallestTerm, GrevLexSmallest, GRevLexLargest, Random, RandomNonzero} each with value an integer (the values need not sum to 100).  If one value is twice the size of another, that strategy will be employed twice as often.  For example, {\tt StrategyDefaultNonRandom} was created by the command:

{{\small\color{blue}
\begin{verbatim}
StrategyDefaultNonRandom = new HashTable from {
    LexLargest => 0,
    LexSmallestTerm => 25,
    LexSmallest => 25,
    GRevLexSmallestTerm => 25,
    GRevLexSmallest => 25,
    GRevLexLargest => 0,
    Random => 0,
    RandomNonzero => 0,
    Points => 0
};
\end{verbatim}}
}

For a tutorial on choosing strategies, see the package documentation, particularly the item:
\begin{center}
    {\tt FastMinorsStrategyTutorial}
\end{center}

\section{Find a submatrix of a given rank: {\tt getSubmatrixOfRank}}
\label{sec.GetSubmatrixOfRank}

This method examines the submatrices of an input matrix and attempts to find one of a 
given rank. If a submatrix with the specified rank is found, a list of two lists is returned. 
The first is the list of row indices, the second is the list of column indices, which describe the desired submatrix of the desired rank. 
If no such submatrix is found, the function will return {\tt null}. 

\par The option {\tt MaxMinors} allows the user to control how many minors to consider before giving up. If left {\tt null}, 
the number considered is based on the size of the matrix. This method will choose the indicated amount of 
minors using one of the strategy options described above. If one of the chosen submatrices has the desired 
rank, the function will terminate and return its rows and columns. This process continues until a submatrix 
is found or {\tt MaxMinors} submatrices have been unsuccessfully checked.  The strategy can be controlled using the {\tt Strategy} option as described above, the default value is {\tt  StrategyDefaultNonRandom}.

\subsection{Examples of {\tt getSubmatrixOfRank}}
In the following example, we first create a $3 \times 4$ matrix over $\bQ[x,y,z]$, $M$. We execute two calls 
to {\tt getSubmatrixOfRank}, the first has no Strategy parameter and the second utilizes 
{\tt StrategyGRevLexSmallest}. Note that these calls return different indices, but both find 
valid rank $3$ submatrices. 
{\small\color{blue}
\begin{verbatim}
  i1 : loadPackage "FastMinors";

  i2 : R = QQ[x,y];
  
  i3 : M = random(R^{2,2,2},R^4)
  
  o3 = {-2} | x2+2/3xy+9/2y2    3/10x2+2/3xy+1/5y2 2x2+5/3xy+7/5y2    4/3x2+1/3xy+10/9y2 |
       {-2} | 3/2x2+2/3xy+2y2   1/2x2+3/2xy+3/4y2  6x2+5xy+4y2        9/5x2+1/5xy+7/2y2  |
       {-2} | 1/4x2+1/7xy+5/6y2 7/5x2+4xy+4/5y2    10/9x2+3/7xy+5/9y2 5/2x2+xy+7/6y2     |
               3       4
  o3 : Matrix R  <--- R
  
  i4 : getSubmatrixOfRank(3,M)
  
  o4 = {{2, 0, 1}, {0, 1, 3}}
  
  o4 : List
  
  i5 : getSubmatrixOfRank(3, M, Strategy=>StrategyGRevLexSmallest)
  
  o5 = {{0, 2, 1}, {1, 2, 0}}
  
  o5 : List
\end{verbatim}
}

In our next example, over a ring with 6 variables, we create a Jacobian matrix out of an ideal generated by 8 random forms of various degrees.
We display the time needed for the {\tt rank} function to return, 
followed by the time elapsed during a call to {\tt getSubmatrixOfRank} when searching for a rank $6$
submatrix.  
We find that {\tt getSubmatrixOfRank} significantly outperformed {\tt rank}. 

{\small\color{blue}
\begin{verbatim}
  i6 : R = ZZ/103[x_1..x_6]

  o6 = R
  
  o6 : PolynomialRing
  
  i7 : J = jacobian ideal apply(8, i -> random(2+random(2), R));
  
               6       8
  o7 : Matrix R  <--- R
  
  i8 : time rank J
       -- used 21.8251 seconds
  
  o8 = 6
  
  i9 : time getSubmatrixOfRank(6, J)
       -- used 0.00714912 seconds
  
  o9 = {{5, 1, 3, 4, 2, 0}, {5, 2, 6, 0, 4, 7}}
  \end{verbatim}
}
In one of the core examples from the {\tt RationalMaps} package, before using this package a function would look at several thousands of submatrices (randomly) typically before finding a submatrix of the desired rank whereas, this package finds one after looking at fewer than half a dozen (typically only looking at 1 or 2 submatrices).  Using this package sped up the computation of that example by more than one order of magnitude, see 
\cite[Page 7]{BottHassanzadehSchwedeSmolkinRationalMaps}, the non-maximal linear rank example.

\section{Finding lower bounds for matrix ranks: {\tt isRankAtLeast}}
\label{sec.IsRankAtLeast}

This method is a direct implementation of {\tt getSubmatrixOfRank}. This function returns a 
boolean value indicating whether the rank of an input matrix, $M$, is greater than or equal to an 
input integer, $n$. In order to do so, the function first performs some basic checks 
to ensure a rank of $n$ is possible given $M$'s dimensions, then executes a call to 
{\tt getSubmatrixOfRank}. If {\tt getSubmatrixOfRank} returns a matrix, then this function 
will return true. However, if {\tt getSubmatrixOfRank} does not return a matrix, a conclusive 
answer can not be reached. As such, the method will then evaluate the rank of $M$ and return 
the appropriate boolean value. 

The function {\tt isRankAtLeast}, which is efficient when {\tt getSubmatrixOfRank} returns quickly, however may be costly if the results are inconclusive and a rank evaluation is necessary. As such, the 
described implementation is not optimized. In order to lead to time improvements, we developed a 
multithreaded version of this function that simultaneously evaluates the rank of $M$ and invokes 
{\tt getSubmatrixOfRank}. Once a thread has terminated with a usable answer, the other threads are cancelled and the
appropriate value is returned.  During the implementation of this functionality, we discovered that 
{\it Macaulay2} becomes unstable when cancelling threads and thus currently do not allow users 
to invoke the multithreaded version. However, this functionality is included in the package and can 
be made easily accessible once the stability issue is resolved. 

\subsection{Example of {\tt isRankAtLeast}}
The following example first creates a $9 \times 9$ 
matrix, $N$, and calls {\tt isRankAtLeast} to determine whether its rank is at least 7. Directly calling 
{\tt rank N} on a matrix of this size would take multiple seconds, whereas {\tt isRankAtLeast} 
returns in a fraction of the time. 

{\small\color{blue}
\begin{verbatim} 
  i1 : loadPackage "FastMinors";

  i2 : N = random(R^{6,6,6,6,6,6,6,7,7},R^9);
                9       9
  o2 : Matrix R  <--- R

  i3 : elapsedTime isRankAtLeast(7,N)
  -- 0.0654172 seconds elapsed

  o3 = true
\end{verbatim}
}

\section{Regular in codimension $n$: {\tt regularInCodimension}}
\label{sec.Rn}

Using the {\tt getSubmatrixOfRank} routines, we provide a function for checking whether a variety is regular in codimension $n$, or $R_n$.  The default strategy is {\tt Strategy=>Default}.

The function {\tt regularInCodimension(ZZ, Ring)} returns {\tt true} if it verifies that the ring is regular in codimension $n$. This only works if the ring is equidimensional, as it is using a Jacobian criterion.  If it cannot make a determination, it returns {\tt null}.  If it ended up computing all minors of the matrix, and it still doesn't have the desired codimension, it will return {\tt false} (note this will likely only occur for small matrices).

\subsection{Example of {\tt regularInCodimension}}
We begin with an example of a 3-dimensional ring that is regular in codimension 1, but not in codimension 2.  It is generated by 12 equations in 7 variables.

{\small\color{blue}
\begin{verbatim}
i3 :  T = ZZ/101[x1,x2,x3,x4,x5,x6,x7];

i4 :  I =  ideal(x5*x6-x4*x7,x1*x6-x2*x7,x5^2-x1*x7,x4*x5-x2*x7,x4^2-x2*x6,x1*x4-x2*x5,
x2*x3^3*x5+3*x2*x3^2*x7+8*x2^2*x5+3*x3*x4*x7-8*x4*x7+x6*x7,x1*x3^3*x5+3*x1*x3^2*x7
+8*x1*x2*x5+3*x3*x5*x7-8*x5*x7+x7^2,x2*x3^3*x4+3*x2*x3^2*x6+8*x2^2*x4+3*x3*x4*x6
-8*x4*x6+x6^2,x2^2*x3^3+3*x2*x3^2*x4+8*x2^3+3*x2*x3*x6-8*x2*x6+x4*x6,x1*x2*x3^3
+3*x2*x3^2*x5+8*x1*x2^2+3*x2*x3*x7-8*x2*x7+x4*x7,x1^2*x3^3+3*x1*x3^2*x5+8*x1^2*x2
+3*x1*x3*x7-8*x1*x7+x5*x7);

o4 : Ideal of T

i5 : S = T/I; dim S

o6 = 3

i7 : time regularInCodimension(1, S)
     -- used 0.150734 seconds

o7 = true

i8 : time regularInCodimension(2, S)
     -- used 2.12777 seconds

i9 : time singularLocus S;
     -- used 8.29746 seconds

i10 : time dim o9
-- used 23.2483 seconds

o10 = 1
\end{verbatim}
}
As seen above, the function {\tt regularInCodimension} verified that $S$ was regular in codimension 1 in a fraction of a second.  When {\tt regularInCodimension(2, S)} was called, nothing was returned, indicating that nothing was found (our function could not make a determination).  Computing the Jacobian ideal however took more than 8 seconds and verifying that it had dimension 1 took more than 23 seconds.  

\subsection{Options and strategies for {\tt regularInCodimension}}

We consider the same example using some different strategies.  For another look at options in this function, see the tutorial in the document under the key:
\begin{center}
    {\tt RegularInCodimensionTutorial}
\end{center}

One might think that it might be just as effective to choose random matrices as to use our strategies, and sometimes it is, but this is not the typical behavior we have observed.
{\scriptsize\color{blue}
\begin{verbatim}
i11 :  time regularInCodimension(1, S, Strategy=>StrategyRandom, Verbose=>true)
regularInCodimension: ring dimension =3, there are 17325 possible minors, we will compute up to 317.599 of them.
regularInCodimension: About to enter loop
internalChooseMinor: Choosing Random
regularInCodimension:  Loop step, about to compute dimension.  Submatrices considered: 7, and computed = 7
regularInCodimension:  isCodimAtLeast failed, computing codim.
regularInCodimension:  partial singular locus dimension computed, = 2
regularInCodimension:  Loop step, about to compute dimension.  Submatrices considered: 9, and computed = 9
regularInCodimension:  isCodimAtLeast failed, computing codim.
regularInCodimension:  partial singular locus dimension computed, = 2
regularInCodimension:  Loop step, about to compute dimension.  Submatrices considered: 11, and computed = 11
regularInCodimension:  isCodimAtLeast failed, computing codim.
regularInCodimension:  partial singular locus dimension computed, = 2
regularInCodimension:  Loop step, about to compute dimension.  Submatrices considered: 14, and computed = 14
regularInCodimension:  isCodimAtLeast failed, computing codim.
regularInCodimension:  partial singular locus dimension computed, = 2
regularInCodimension:  Loop step, about to compute dimension.  Submatrices considered: 18, and computed = 18
regularInCodimension:  isCodimAtLeast failed, computing codim.
regularInCodimension:  partial singular locus dimension computed, = 2
regularInCodimension:  Loop step, about to compute dimension.  Submatrices considered: 24, and computed = 24
regularInCodimension:  isCodimAtLeast failed, computing codim.
regularInCodimension:  partial singular locus dimension computed, = 2
regularInCodimension:  Loop step, about to compute dimension.  Submatrices considered: 31, and computed = 31
regularInCodimension:  isCodimAtLeast failed, computing codim.
regularInCodimension:  partial singular locus dimension computed, = 2
regularInCodimension:  Loop step, about to compute dimension.  Submatrices considered: 40, and computed = 40
regularInCodimension:  isCodimAtLeast failed, computing codim.
regularInCodimension:  partial singular locus dimension computed, = 2
regularInCodimension:  Loop step, about to compute dimension.  Submatrices considered: 52, and computed = 52
regularInCodimension:  isCodimAtLeast failed, computing codim.
regularInCodimension:  partial singular locus dimension computed, = 2
regularInCodimension:  Loop step, about to compute dimension.  Submatrices considered: 67, and computed = 67
regularInCodimension:  isCodimAtLeast failed, computing codim.
regularInCodimension:  partial singular locus dimension computed, = 2
regularInCodimension:  Loop step, about to compute dimension.  Submatrices considered: 87, and computed = 87
regularInCodimension:  isCodimAtLeast failed, computing codim.
regularInCodimension:  partial singular locus dimension computed, = 1
regularInCodimension:  Loop completed, submatrices considered = 87, and computed = 87.  
singular locus dimension appears to be = 1
     -- used 1.04945 seconds

o11 = true
\end{verbatim}
}
\noindent
Above, we have deleted 86 of the 87 times the verbose output displays {\tt internalChooseMinor: Choosing Random}.

In this particular example, the {\tt StrategyRandom} option looked at 87 submatrices of the Jacobian matrix.  Note it does not check to see whether we have obtained the desired codimension after considering each new random submatrix.  Instead, it only computes the codimension periodically, with the space between checks increasing.  The {\tt considered} values on each line tells how many submatrices have been considered.  The {\tt computed} tells how many were not repeats ({\tt computed} will be nearly the same as {\tt considered} with a random strategy).

Running {\tt Rn(1, S, Strategy=>StrategyRandom, Verbose=>true)} 50 times yielded:
\begin{enumerate}
  \item 61.3 average number of submatrices of the Jacobian matrix considered.
  \item a median value of 40 or 52 submatrices of the Jacobian matrix considered.
  \item a minimum value of 7 submatrices of the Jacobian matrix considered (one time).
  \item a maximum value of 248 submatrices of the Jacobian matrix considered (one time).
\end{enumerate}
Because of certain settings, we will not check the codimension of the singular locus until 7 submatrices have been considered.  Users can control this behavior via the {\tt MinMinorsFunction} and {\tt CodimCheckFunction} options, see the tutorial in the documentation.

On the other hand, the default strategy {\tt Rn(1, S, Strategy=>StrategyDefaultNonRandom, Verbose=>true)} run 50 times yields
\begin{enumerate}
  \item 12.1 average number of submatrices of the Jacobian matrix considered.
  \item a median value of 7 or 9 submatrices of the Jacobian matrix considered.
  \item a minimum value of 7 submatrices of the Jacobian matrix considered (25 times).
  \item a maximum value of 40 submatrices of the Jacobian matrix considered (one time).
\end{enumerate}

In the above example, {\tt Strategy=>StrategyLexSmallest} yields even better performance.  

Finally, using {\tt Strategy=>StrategyPoints} (combined with the options {\tt MinMinorsFunction=>(t->t)} and {\tt CodimCheckFunction=>(t->t)}) to check codimension after computing every submatrix, produces:
\begin{enumerate}
  \item 4.96 average number of submatrices of the Jacobian matrix considered.
  \item a median value of 5 submatrices of the Jacobian matrix considered.
  \item a minimum value of 4 submatrices of the Jacobian matrix considered (3 times).
  \item a maximum value of 6 submatrices of the Jacobian matrix considered (one time).
\end{enumerate}
In this case, {\tt StrategyPoints} considers very few submatrices, but it actually does the computation substantially slower than {\tt StrategyDefaultNonRandom} since finding each submatrix can be a lot of work as rational points must be found.  However, {\tt StrategyPoints} is still faster than {\tt StrategyRandom}.

Note that larger matrices tend to exhibit even larger disparities between the strategies.

\subsection{Notes on implementation}

As mentioned above, this function computes minors (based on the passed {\tt Strategy}) option until either it finds that the singular locus has the desired dimension, or until it has considered too many minors.  By default, it considers up to:
\[
  10\cdot (\text{Minimum number of minors needed}) + 8 \cdot \log_{1.3}(\text{possible minors}).
\]
This was simply chosen by experimentation.
If the user is trying to show a singular locus has a certain codimension, they will need a minimum number of minors.  The multiplication by 10 is due to our default strategy using multiple strategies, but only considering one might work well on a given matrix.  The user can set the option {\tt MaxMinors} to a function $F$ with two inputs, $x = (\text{minors needed})$ and $y = (\text{possible minors})$, where $F$ outputs the maximum number of matrices to compute.  More simply, one may simply set {\tt MaxMinors} a number.

These matrices are considered in a loop.  We begin with computing a constant number of minors, by default $2 \cdot (\text{Minimum number of minors needed}) + 3$, and check whether the output has the right dimension.  The user can also set the option {\tt MinMinorsFunction} to a function $G$ with one input, $x = (\text{minors needed})$, which will output how many minors to compute before first checking the codimension.
After those initial minors are found, we compute additional minors, checking periodically (based on an exponential function, $1.3^k$ minors considered before the next reset) whether our minors define a subset of the desired codimension.  New functions can be provided via the option {\tt CodimCheckFunction}, see the tutorial for more details.  If in this loop, a submatrix is considered again, it is not recomputed, but the counter is still increased.

\subsection{Other options}

This function also includes other options including the option {\tt Modulus} which handles switching the coefficient field for a field of characteristic $p > 0$ (which is specified with {\tt Modulus => p}. )

One can also control how determinants are computed with the {\tt DetStrategy} option, valid values are {\tt Bareiss}, {\tt Cofactor} and {\tt Recursive}.

\section{Projective dimension: {\tt projDim}}
\label{sec.ProjDim}

In April of 2019, it was pointed out in a thread on github
\begin{center}
  https://github.com/Macaulay2/M2/issues/936
\end{center}
that the command {\tt pdim} sometimes provides an incorrect value (an overestimate) for projective dimension for non-homogeneous modules over polynomial rings.  There it was also suggested that this could be addressed by looking at appropriate minors of the matrices in a possibly non-minimal resolution, but that in practice these matrices have too many minors to compute.  We have implemented a function {\tt projDim} that tries to address this by looking at only \emph{some} minors.  Our function does not solve the problem as it also gives only an upper bound on the projective dimension.  However, this upper bound is more frequently correct.

The idea is as follows.  Take a free resolution of a module $M$ over a polynomial ring $R$  
\[
  \xymatrix{
  0 & \ar[l] M  & \ar[l] F_0 & \ar[l]_{d_1} F_1 & \ar[l] \dots & \ar[l]_{d_{n-1}} F_{n-1}  & \ar[l]_{d_n} F_n & \ar[l] 0.  
  }
\] 
Each $d_i$ is given by a matrix.  The term $F_n$ is unnecessary (i.e., $d_n$ splits) happens exactly when the $\rank F_n$ minors of $d_n$ generate the unit ideal.  In that case, we know or projective dimension is at most $n-1$.  However, we can continue in this way, we can compute the $(\rank F_{n-1} - \rank F_n)$-minors of $d_{n-1}$, and see whether they generate the unit ideal.  Our algorithm of course only computes a subset of those minors.

\subsection{Example of {\tt projDim}}
In the below example, we take a monomial ideal of projective dimension 2, compute a non-homogeneous change of coordinates, and observe that {\tt pdim} returns an incorrect answer that {\tt projDim} corrects.

{\small\color{blue}
\begin{verbatim}
i1 : R = QQ[x,y,z,w];

i2 : I = ideal(x^4,x*y,w^3, y^4);

i3 : pdim module I

o3 = 2

i4 : f = map(R, R, {x+x^2+1, x+y+1, z+z^4+x-2, w+w^5+y+1});

i5 : pdim module f I

o5 = 3

i6 : time projDim module f I
    -- used 3.43851 seconds

o6 = 2

i7 : time projDim(module f I, MinDimension=>2)
    -- used 0.0503165 seconds

o7 = 2
\end{verbatim}
}

\subsection{Options}

As seen in the previous example, setting {\tt MinDimension} will can substantially speed up the computation as otherwise, the function will try to determine whether the projective dimension is actually $1$.  

The option {\tt MaxMinors} can either be set to the number of a minors computed
      at each step.  Alternatively, it can be set to be a list of numbers, one
      for each step in the above algorithm. Finally, it can be set to be a function of the dimension $d$ of the polynomial ring $R$ and the number $t$ of possible minors.  This is the default option, and the function is:  $5*d + 2*\log_{1.3}(t)$.
The option {\tt Strategy} is also available and it works as in the above functions with the default value being {\tt StrategyDefault}.

\section{Computing ideals of minors: {\tt recursiveMinors}}
\label{sec.RecursiveMinors}
{\it Macaulay2} contains a {\tt minors} method that returns the ideal 
of minors of a certain size, $n$,  in a given matrix, a necessary step in 
locating singularities. However, the current implementation’s default is to 
evaluate determinants using the Bareiss algorithm, which is efficient when 
the entries in the matrix have a low degree and few variables, but very slow 
otherwise. The current minors method also allows users to compute determinants 
using cofactor expansion, but this strategy performs some unnecessary calculations, 
causing it to be quite costly as well. We improved the current cofactor 
expansion method to find the determinants of minors by adding recursion and 
multithreading throughout. We also eliminated said unnecessary calculations by 
ensuring that only the required determinants are being computed at each step 
of the recursion, rather than all possible determinants of the given size.

In order to do so, we programmed a method in {\it Macaulay2}’s software that 
recursively finds all $n \times n$ minors by first computing the $2 \times 2$ minors 
and storing them in a hash table. Then we use the $2 \times 2$ minors to compute the 
necessary $3 \times 3$ minors, and so forth. This process is repeated recursively until 
the minors of size $n \times n$ are evaluated. At each step, we only compute the 
determinants that will be needed when performing a cofactor expansion on the 
following size minor.

To allow for further time improvements, we also utilized {\it Macaulay2}’s existing 
parallel programming methods to multithread our code so different computations 
at each step of the recursion can occur simultaneously in separate threads. We 
divide the list of all determinants to be evaluated into different available 
threads and wait for them to finish before consolidating the results in a hash 
table and proceeding with the recursion. In order to more effectively utilize 
{\it Macaulay2}'s multithreading methods, we also created a nanosleep method that 
waits a given number of nanoseconds, rather than full seconds. This function has already 
been incorporated into the software. 

\subsection{Example of {\tt recursiveMinors}}
\par Below, we first create a simple matrix, M, of polynomials in a single variable with 
rational coefficients and execute the {\tt recursiveMinors} method to find the ideal of all 
$3 \times 3$ minors. As can be seen, the result is equivalent to the output of the minors method 
when called with the same parameters. We then create a new, larger matrix, $N$, with two
dimensional rational coefficients and return the computation time for {\tt recursiveMinors} and 
minors utilizing both the Bareiss and Cofactor strategies. The {\tt recursiveMinors} method 
finished executing approximately six times faster than the Bareiss algorithm and almost seven
times faster than the Cofactor expansion, while yielding the same results. 

{\small\color{blue}
\begin{verbatim}
  i1 : loadPackage "FastMinors";  

  i2 : allowableThreads => 8;

  i3 : R = QQ[x];

  i4 : M = random(R^{2,2,2}, R^4)

  o4 = {-2} | x2    3x2   5/8x2 7/10x2 |
       {-2} | 3/4x2 2x2   7/4x2 9x2    |
       {-2} | x2    2/9x2 1/2x2 4/3x2  |
               3       4
  o4 : Matrix R  <--- R
  
  i5 : recursiveMinors(3,M)

              1403 6  449 6    292 6  517 6
  o5 = ideal (----x , ---x , - ---x , ---x )
               60     240       45    144
  
  o5 : Ideal of R
   
  i6 : recursiveMinors(3,M) == minors(3,M)
  
  o6 = true
  
  i7 : Q = QQ[x,y];

  i8 : N = random(Q^{5,5,5,5,5,5}, Q^7);
                6       7
  o8 : Matrix Q  <--- Q

  i9 : elapsedTime minors(5,N, Strategy => Bareiss);
  -- 1.42867 seconds elapsed
  
  o9 : Ideal of Q
  
  i10 : elapsedTime minors(5,N, Strategy => Cofactor);
  -- 1.82251 seconds elapsed
  
  o10 : Ideal of Q
  
  i11 : elapsedTime recursiveMinors(5,N);
  -- 0.273007 seconds elapsed;
  
  o11 : Ideal of Q
  
  i12 : recursiveMinors(5,N) == minors(5,N)
  
  o12 = true

\end{verbatim}
}

We briefly give a table showing the limits of this package.  We consider a random $6 \times 7$ matrix over $\bQ[x,y]$ as above, and then also for $\bQ[x,y,z]$.  We compare the single-threaded and 4-threaded version of {\tt recursiveMinors} in this package with the {\tt Bareiss} and {\tt Cofactor} strategies with {\tt recursiveMinors} for different degrees of the terms.

\begin{figure}[t]
    \begin{tabular}{c||c|c|c|c}
        {\bf Degree} & {\bf Bareiss} & {\bf Cofactor} & {\bf RecursiveMinors} & {\bf RecursiveMinors, Threads$=>$4} \\ \hline \hline
        {\tt 8} & 3.46466 sec & 4.44334 sec & 0.631749 sec & 0.407553 sec \\
        {\tt 10} & 5.77108 sec & 6.79853 sec & 0.971093 sec & 0.55966 sec \\
        {\tt 12} & 7.40464 sec & 8.93534 sec & 1.21964 sec & 0.699108 sec \\
        {\tt 15} & 12.1871 sec &  12.0075 sec &  1.68724 sec & 1.00743 sec \\
        {\tt 20} & 21.0065 sec &  22.6152 sec &  2.81906 sec & 1.85354 sec \\
        {\tt 25} & 31.915 sec &  34.8651 sec &  4.23277 sec & 2.63514 sec \\
        {\tt 40} & 83.5833 sec & 77.1983 sec &   10.5852 sec & 6.29589 sec \\
        {\tt 60} & 181.179 sec &  192.911 sec &   23.875 sec & 13.0619 sec \\
    \end{tabular}
    \caption{Time to compute the $5 \times 5$ minors of a $6 \times 7$ random matrix over $\bQ[x,y]$ where the entries have degree {\bf Degree}.  }
\end{figure}

Generally speaking, {\tt reursiveMinors} performs best when the matrix one is looking at has very expensive-to-compute minors (such as with the random matrices we consider above).  In sparse examples, and examples with easy-to-compute determinants, other strategies tend to perform better.

\begin{figure}[h]
    \begin{tabular}{c||c|c|c|c}
        {\bf Degree} & {\bf Bareiss} & {\bf Cofactor} & {\bf RecursiveMinors} & {\bf RecursiveMinors, Threads$=>$4} \\ \hline \hline
        {\tt 2} & 5.99812 sec & 3.78488 sec & 0.588191 sec & 0.51851 sec \\        
        {\tt 3} & 17.3966 sec & 8.78131 sec & 1.73022 sec & 1.53543 sec \\        
        {\tt 4} & 49.6152 sec & 22.5747 sec & 4.58156 sec &  3.83319 sec \\  
        {\tt 5} & 115.412 sec & 45.0877 sec & 8.36395 sec &  6.39432 sec \\  
    \end{tabular}
    \caption{Time to compute the $5 \times 5$ minors of a $6 \times 7$ random matrix over $\bQ[x,y,z]$ where the entries have degree {\bf Degree}.  }
\end{figure}

\section{Performance and limits of the package}
\label{sec.Performance}

We conclude by providing some figures showing how long various computations take in several different strategies.  We limit ourselves to the function {\tt regularInCodimension} as other functions such as {\tt projDim} have roughly similar performance.  Note some discussion of the performance behavior of taking determinants (including via a recursive algorithm) was already discussed above.  Again, we recommend the interested user also see the tutorial 
\begin{center}{\tt FastMinorsStrategyTutorial}\end{center} in the package documentation.

The column {\bf Attempts} indicates how many time we ran this computation.  The column {\bf Successful} shows what percentage of the time the function verified that the given equation was regular in a certain codimension (depending on the strategy, it doesn't always succeed).  All computations were run in {\it Macaulay2} version 1.18 on a machine running Ubuntu 20.04 with 64 gigabytes of memory.

In \autoref{fig.R1ConeOverProductOfElliptic}, we verify that the cone over a product of elliptic curves (an Abelian surface) embedded in $\bP^8$ is regular in codimension 1.   Note that {\tt StrategyRandom} does not tend to work well on this or other examples, and so we generally do not consider it further.  In \autoref{fig.R2ConeOverProductOfElliptic} we verify the same example is regular in codimension 2. When we make the elliptic curves defined by less sparse equations, {\tt Points} tends to perform much better, as can be seen in in \autoref{fig.R1ConeOverProductOfEllipticNonSparse}. 

We next consider a relatively sparse higher dimension example in \autoref{fig.R1ConeOverProductOfEllipticAndP1}.  Here we are taking a cone over a product of an elliptic curve with a diagonal equation, an elliptic curve in Weierstrass form and a copy of $\bP^1$.  This is a cone over a 3-dimensional smooth projective variety embedded in $\bP^{17}$.

\begin{figure}[h]
    \begin{tabular}{c||c|c|c}
        {\bf Strategy} & {\bf Attempts} & {\bf Average time} & {\bf Successful}\\\hline \hline
        {\tt StrategyDefault} & 100 & 1.7 sec & 100\% \\\hline
        {\tt StrategyDefaultNonRandom} & 100 & 0.9 sec & 100\% \\\hline
        {\tt Points} & 100 & 4.0 sec & 100\%\\\hline
        {\tt StrategyDefaultWithPoints} & 100 & 2.2 sec & 100\%\\\hline
        {\tt StrategyRandom} & 100 & 6.1 sec & {\it 4\%} \\\hline
        {\tt StrategyRandom, MaxMinors=>2000} & 20 & 49.0 sec & {\it 15\%} \\\hline        
        {\tt StrategyRandom, MaxMinors=>5000} & 10 & 238.1 sec & {\it 50\%} \\\hline        
    \end{tabular}
    \begin{center}Regular in codimension 1, 9 variables, 28 equations, 31,646,160 possible 6 by 6 minors\end{center}
    \caption{We check $R$ is regular in codimension 1 where $R$ is the cone over a product of two elliptic curves in positive characteristic given with a Segre embedding.  One of the curves is diagonal, the other is in Weierstrass form.  This has a relatively sparse Jacobian matrix.}
    \label{fig.R1ConeOverProductOfElliptic}
\end{figure}

\begin{figure}[h]
    \begin{tabular}{c||c|c|c}
        {\bf Strategy} & {\bf Attempts} & {\bf Average time} & {\bf Successful}\\\hline \hline
        {\tt StrategyDefault} & 10 & 10.9 sec & {\it 0\%} \\\hline
        {\tt StrategyDefault, MaxMinors=>5000} & 10 & 30.1 sec & {100\%} \\\hline
        {\tt StrategyDefaultNonRandom} & 10 & 7.7 sec & {\it 0\%} \\\hline
        {\tt StrategyDefaultNonRandom, MaxMinors=>5000} & 10 & 13.7 sec & 100\% \\\hline
        {\tt Points} & 10 & 4.6 sec & 100\%\\\hline
        {\tt StrategyDefaultWithPoints} & 10 & 5.8 sec & 100\%\\\hline
    \end{tabular}
    \begin{center}Regular in codimension 2, 9 variables, 28 equations, 31,646,160 possible $6 \times 6$ minors\end{center}
    \caption{We check $R$ is regular in codimension 2 where $R$ is the cone over a product of two elliptic curves in positive characteristic given with a Segre embedding.  One of the curves is diagonal, the other is in Weierstrass form.  This has a relatively sparse Jacobian matrix.  Using {\tt StrategyDefault} and {\tt StrategyDefaultNonRandom} did not work with the default number of minors, but increasing {\tt MaxMinors} led to successful verification that the ring was regular in codimension 2.}
    \label{fig.R2ConeOverProductOfElliptic}
\end{figure}

\begin{figure}[h]
    \begin{tabular}{c||c|c|c}
        {\bf Strategy} & {\bf Attempts} & {\bf Average time} & {\bf Successful}\\\hline \hline
        {\tt StrategyDefault} & 10 & $\infty?$ sec & {\it 0\%} \\\hline
        {\tt StrategyDefaultNonRandom} & 10 & $\infty?$ sec & {\it $<$10\%} \\\hline
        {\tt Points, CodimCheckFunction => t->t+1} & 10 & 7.7 sec & 100\%\\\hline
        {\tt StrategyDefaultWithPoints} & 10 & ? sec & about 50\%\\\hline
    \end{tabular}
    \begin{center}Regular in codimension 1, 9 variables, 28 equations, 31,646,160 possible $6 \times 6$ minors\end{center}
    \caption{We check $R$ is regular in codimension 1 where $R$ is the cone over a product of two elliptic curves in positive characteristic given with a Segre embedding.  One of the curves is in Weierstrass form, the other is given by a random degree 3 equation.  This has a relatively complicated (non-sparse) Jacobian matrix.  The other strategies generally do not work.  The one exception is {\tt StrategyDefaultWithPoints} which sometimes is very fast (faster than {\tt Points}), and other times gets stuck trying to compute a point.  Setting {\tt CodimCheckFunction => t->t+1} forces the codimension to be checked at every step, which provides better and more consistant performance.  Without that, sometimes this function will hang trying to find a point after on a 1-dimensional scheme where it has already verified that $R$ is regular in codimension 1, but has not computed that codimension yet.}
    \label{fig.R1ConeOverProductOfEllipticNonSparse}
\end{figure}

\begin{figure}[h]
    \begin{tabular}{c||c|c|c}
        {\bf Strategy} & {\bf Attempts} & {\bf Average time} & {\bf Successful}\\\hline \hline        
        {\tt StrategyDefaultNonRandom} & 10 & $\infty?$ sec & {\it 0\%} \\\hline
        {\tt Points} & 10 & 58.5 sec & 100\%\\\hline
        {\tt StrategyDefaultWithPoints} & 10 & 27.13 sec & 100\%\\\hline
    \end{tabular}
    \begin{center}Regular in codimension 1, 18 variables, 139 equations,17,927,476,818,965,522,386,560 possible $14 \times 14$ minors\end{center}
    \caption{We check $R$ is regular in codimension 1 where $R$ is the cone over a product of two elliptic curves plus a $\bP^1$ in positive characteristic given with a Segre embedding.  One of the curves is in Weierstrass form, the other is given by a random degree 3 equation.  This has a relatively sparse Jacobian matrix.  The other strategies (not involving points) generally do not work.  Using {\tt StrategyDefaultNonRandom} took more than 30 minutes and computed more than 5000 minors, but still did not finish.}
    \label{fig.R1ConeOverProductOfEllipticAndP1}
\end{figure}

We now move on to computing dimensions of singular loci of varieties that are not cones.  We constructed several non-normal (non-S2) varieties using the {\tt Pullback} package \cite{PullbackSource}.  First, in \autoref{fig.R1Pullback3LinesGluedInA3} we took 3 coordinate axes through the origin in $\bA^3$ and randomly glued them to a single line.  In \autoref{fig.R1Pullback3RandomLinesGluedInA3} we did the same with three random lines through the origin (creating a less sparse Jacobian matrix).  Finally, in \autoref{fig.R1Pullback3LinesGluedInA4}, we consider a similar example in $\bA^4$ (except now it is regular in codimension 2), first verifying it is regular in codimension 1.  Finally, we verify it is regular in codimension 2 in \autoref{fig.R2Pullback3LinesGluedInA4} .  

\begin{figure}[h]
    \begin{tabular}{c||c|c|c}
        {\bf Strategy} & {\bf Attempts} & {\bf Average time} & {\bf Successful}\\\hline \hline        
        {\tt StrategyDefault} & 100 & 0.8 sec & 100\% \\\hline
        {\tt StrategyDefaultNonRandom} & 100 & 0.5 sec & 100\% \\\hline
        {\tt Points} & 100 & 3.0 sec & 100\%\\\hline
        {\tt StrategyDefaultWithPoints} & 100 & 2.4 sec & 100\%\\\hline
    \end{tabular}
    \begin{center}Regular in codimension 1, 8 variables, 26 equations, 3,683,680 possible $5 \times 5$ minors\end{center}
    \caption{We check $R$ is regular in codimension 1 where $R$ is obtained by gluing three coordinate axis lines through the origin in $\bA^3$ together to a single line.  This is a 3-dimensional ring that is regular in codimension 1, but not codimension 2.  The Jacobian matrix is fairly sparse, but has some quite complicated sections.}
    \label{fig.R1Pullback3LinesGluedInA3}
\end{figure}

\begin{figure}[h]
    \begin{tabular}{c||c|c|c}
        {\bf Strategy} & {\bf Attempts} & {\bf Average time} & {\bf Successful}\\\hline \hline        
        {\tt StrategyDefault} & 20 & 2.0 sec & 100\% \\\hline
        {\tt StrategyDefaultNonRandom} & 20 & 0.5 sec & 100\% \\\hline
        {\tt Points} & 10 & 11.6 sec & 100\%\\\hline
        {\tt StrategyDefaultWithPoints} & 10 & 6.9 sec & 100\%\\\hline
    \end{tabular}
    \begin{center}Regular in codimension 1, 8 variables, 34 equations, 15,582,336 possible $5 \times 5$ minors\end{center}
    \caption{We check $R$ is regular in codimension 1 where $R$ is obtained by gluing random lines through the origin in $\bA^3$ together to a single line.  This is a 3-dimensional ring that is regular in codimension 1, but not codimension 2.  The Jacobian matrix is substantially less sparse than when we glued the three \emph{coordinate axes}.}
    \label{fig.R1Pullback3RandomLinesGluedInA3}
\end{figure}

\begin{figure}[h]
    \begin{tabular}{c||c|c|c}
        {\bf Strategy} & {\bf Attempts} & {\bf Average time} & {\bf Successful}\\\hline \hline        
        {\tt StrategyDefault} & 100 & 5.2  sec & 100\% \\\hline
        {\tt StrategyDefaultNonRandom} & 100 & 1.3 sec & 100\% \\\hline
        {\tt Points} & 20 & 7.3 sec & 100\%\\\hline
        {\tt StrategyDefaultWithPoints} & 20 & 4.0 sec & 100\%\\\hline
    \end{tabular}
    \begin{center}Regular in codimension 1, 11 variables, 52 equations, 44,148,904,800 possible $7 \times 7$ minors\end{center}
    \caption{We check $R$ is regular in codimension 1 where $R$ is obtained by gluing three coordinate axis lines through the origin in $\bA^4$ together to a single line.  This is a 4-dimensional ring that is regular in codimension 2, but not codimension 3.  The Jacobian matrix is fairly sparse, but has some quite complicated sections.}
    \label{fig.R1Pullback3LinesGluedInA4}
\end{figure}

\begin{figure}[h]
    \begin{tabular}{c||c|c|c}
        {\bf Strategy} & {\bf Attempts} & {\bf Average time} & {\bf Successful}\\\hline \hline        
        {\tt StrategyDefault} & 20 & 14.9  sec & 100\% \\\hline
        {\tt StrategyDefaultNonRandom} & 20 & 5.5 sec & 100\% \\\hline
        {\tt Points} & 10 & $\infty?$ sec & {\it 0\%}\\\hline
        {\tt StrategyDefaultWithPoints} & 10 & $\infty?$ sec & {\it 0\%}\\\hline        
    \end{tabular}
    \begin{center}Regular in codimension 2, 11 variables, 52 equations, 44,148,904,800 possible $7 \times 7$ minors\end{center}
    \caption{We check $R$ is regular in codimension 2 where $R$ is obtained by gluing three coordinate axis lines through the origin in $\bA^4$ together to a single line.  This is a 4-dimensional ring that is regular in codimension 2, but not codimension 3.  The Jacobian matrix is fairly sparse, but has some quite complicated sections.  Strategies involving {\tt Points} fail quickly as they use more than 64 gigabytes of ram.}
    \label{fig.R2Pullback3LinesGluedInA4}
\end{figure}

 \clearpage
\bibliographystyle{skalpha}
\bibliography{MainBib}

\end{document}